\documentclass[11pt,reqno]{amsart}
\usepackage{graphicx, amssymb, color,pdfsync}
\usepackage[all,cmtip]{xy}
\numberwithin{equation}{section}
\usepackage{mathrsfs}
\usepackage{hyperref}
\usepackage{diagbox}
\usepackage{booktabs}
\usepackage{multirow}
\usepackage{rotating}
\usepackage{tikz}
\usetikzlibrary{matrix,arrows}
\DeclareMathOperator{\Aut}{Aut}
\usepackage[switch]{lineno}

\begin{document}
\newcommand{\s}{\vspace{0.2cm}}

\newtheorem{theo}{Theorem}
\newtheorem{prop}{Proposition}
\newtheorem{coro}{Corollary}
\newtheorem{lemm}{Lemma}
\newtheorem{claim}{Claim}
\newtheorem{example}{Example}
\theoremstyle{remark}
\newtheorem{rema}{\bf Remark}
\newtheorem*{rema-non}{\bf Remark}
\newtheorem{defi}{\bf Definition}

\title[Riemann surfaces of genus $1+q^2$ with $3q^2$ automorphisms]{Riemann surfaces of genus $1+q^2$ \\with $3q^2$ automorphisms}
\date{}

\author{Angel Carocca and Sebasti\'an Reyes-Carocca}
\address{Departamento de Matem\'atica y Estad\'istica, Universidad de La Frontera, Avenida Francisco Salazar 01145, Temuco, Chile.}
\email{angel.carocca@ufrontera.cl, sebastian.reyes@ufrontera.cl}

\thanks{Partially supported by Fondecyt Grants 1200608, 11180024, 1190001 and Redes Grant 170071}

\begin{abstract} 
In this article we classify compact Riemann surfaces of genus $1+q^2$ with a group of automorphisms of order $3q^2,$ where $q$ is a prime number. We also study decompositions of the corresponding Jacobian varieties.
\end{abstract}
\maketitle

\section{Introduction and statement of the results}

The study and classification of  groups of automorphisms of compact Riemann surfaces (or complex projective algebraic curves) is a classical  problem which has attracted considerable interest since a long time. Regarding this issue, in the late nineteenth century a fundamental result was obtained by Hurwitz, who succeeded in proving that the full automorphism group of a compact Riemann surface of genus $g \geqslant 2$ is finite, and that its order is at most $84g-84.$ Later, this problem acquired a new relevance when its relationship with Teichm\"{u}ller and moduli  spaces was developed.

\s

An interesting problem is to study and describe those compact Riemann surfaces whose automorphism groups share a common property.  The most prominent example concerning that are Hurwitz curves; namely, those Riemann surfaces possessing the maximal possible number of automorphisms. Nowadays, it is classically known that Hurwitz curves correspond to  regular covers of the projective line ramified over three values, marked with 2, 3 and 7. Another well-known  example is the cyclic case, which was considered, among others, by Wiman. In \cite{Wi}, he showed that the largest cyclic group of automorphisms of a  Riemann surface of genus $g \geqslant 2$ has order at most $4g+2.$ Furthermore,  the so-called  Wiman curve of type I given by\begin{equation*} \label{eWiman}y^2=x^{2g+1}-1\end{equation*}shows that this upper bound is attained for each $g$; see \cite{Harvey1}. In the early nineties, the uniqueness problem was addressed by  Kulkarni who proved in \cite{K1} that, for $g$ sufficiently large, the aforementioned curve is the unique Riemann surface of genus $g$ with an automorphism of order $4g+2.$

\s

Let $a,b \in \mathbb{Z}.$ Following \cite{Kis}, the sequence $ag+b$ for $g=2,3, \ldots$ is called admissible if for infinitely many values of $g$ there is a compact Riemann surface of genus $g$ with a group of automorphisms of order $ag+b.$ 

In addition to the already mentioned admissible sequences $84g-84$ and $4g+2,$ the classical case $8g+8$ was considered by Accola \cite{Accola},  Maclachlan \cite{Mac} and Kulkarni \cite{K1}. Very recently, the cases $4g+4$ and $4g$ have been studied by Bujalance, Costa and Izquierdo in \cite{BCI} and by Costa and Izquierdo in \cite{CI} respectively; see also \cite{yojpaa}.

\s

Let $\lambda$ be a positive  integer. Belolipetsky and Jones in \cite{BJ} studied the admissible sequence $\lambda(g-1)$ and succeeded in proving that, under the assumption that $g-1$ is a sufficiently large prime number and $\lambda \geqslant 7$, a compact Riemann surface of genus $g$ with a group of automorphisms of order $\lambda(g-1)$ lies in one of six infinite well-described sequences of examples (similar results but stated in a combinatorial point of view of regular maps can be found in \cite{pp}). Later, under the assumption that $g-1$ is prime, the case $\lambda=4$ was classified in \cite{yoisrael} and the cases $\lambda=3,5,6$ were classified in \cite{IRC}; see also  \cite{IJRC} for a unified treatment of this case.

\s

In this article we shall deal with the admissible sequence $3(g-1)$  and extend the results proved in \cite{IRC} by providing a complete classification in the case that $g-1$ is assumed to be the square of a prime number, instead of a prime number. In other words, we study and classify all those compact Riemann surfaces of genus $1+q^2$ endowed with a group of automorphisms of order $3q^2,$ where $q$ is a prime number.

\s

The main result of this paper is the following:

\begin{theo} \label{t1}
Let $q \geqslant 7$ be a prime, set $g=1+q^2$ and let $S$ be a compact Riemann surface of genus $g$ with a group of automorphisms of order $3q^2.$ 

\s

If $q \equiv -1 \mbox{ mod } 3$ then $S$ belongs to  
the complex one-dimensional family $\mathscr{C}_g$ of   compact  Riemann surfaces of genus $g$ with a group of automorphisms isomorphic to $$G_1=\langle a,b, t: a^q=b^q=t^3=[a,b]=1, tat^{-1}=b, tbt^{-1}=(ab)^{-1} \rangle,$$acting with signature $(0;3,3,3,3).$

\s

If $q \equiv 1 \mbox{ mod } 3$ then $S$ belongs to  either:
\begin{enumerate}
\s

\item the family $\mathscr{C}_g$ 
\s

\item the complex one-dimensional family $\mathscr{U}_g$ of compact Riemann surfaces of genus $g$ with a group of automorphisms isomorphic to $$G_2=\langle a, t : a^{q^2}=t^3=1, tat^{-1}=a^{s}\rangle$$where $s$ is a primitive third root of unity in $\mathbb{Z}_{q^2}$, acting with signature $(0; 3,3,3,3),$ or 

\s

\item the complex one-dimensional family $\mathscr{V}_g$ of compact Riemann surfaces of genus $g$ with a group of automorphisms isomorphic to $$G_3=\langle a, b, t : a^{q}=b^q=t^3=[a,b]=1, tat^{-1}=a^{r}, tbt^{-1}=b^r\rangle$$where $r$ is a primitive third root of unity in $\mathbb{Z}_q$, acting with signature $(0; 3,3,3,3).$
\end{enumerate}
\end{theo}

It is worth mentioning that for $q=5$ the family $\mathscr{C}_{26}$ also exists. As a matter of fact, if $S$ is a compact Riemann surface of genus 26 with a group of automorphisms of order $75$ then either $S$ belongs to $\mathscr{C}_{26}$ or admits a triangle action of $C_5 \times C_{15}$ with signature $(0; 5,15,15)$; see \cite{C}.    By contrast, as the reader could expect, the behavior for $q=2$ and $3$ is completely different. For instance, it is straghtforward to verify the existence of  a one-dimensional family of Riemann surfaces of genus ten with action of $C_9 \rtimes C_3$ with signature $(1;3)$ and the existence of a one-dimensional family of Riemann surfaces of genus five with action of $\mathbf{D}_{6} 
$ with signature $(0; 2,2,6,6).$  

\s

The following corollary is a direct consequence of Theorem \ref{t1} together with the fact that each group of order $6q^2$ has  a subgroup  of order $3q^2,$ for $q \geqslant 7$ prime.

\begin{coro}

Let $q \geqslant 7$ be a prime and set $g=1+q^2.$ If $S$ is a compact Riemann surface of genus $g$ with a group of automorphisms of order $6q^2$ then $S$ belongs to one of the families introduced in Theorem \ref{t1}.
\end{coro}

The following three results describe the full automorphism group of the Riemann surfaces lying in the before introduced families, and also the families themselves  (seen as subvarieties of the moduli space) in terms of the number of equisymmetric strata.

\begin{theo} \label{t21}
Let $q \geqslant 5$ be  prime  and set $g=1+q^2.$ The family $\mathscr{C}_g$  consists of at most $(q^2+1)(q+1)$  strata. Moreover, up to finitely many exceptions, the full automorphism group of $S \in \mathscr{C}_g$ is either: \begin{enumerate}
\item isomorphic to $G_1,$ 
\item isomorphic to the group of order $6q^2$  presented as  $$H_1 =\langle a,b,t,z : a^q=b^q=t^3=z^2=[a,b]=[z,t]=1,$$ $$ tat^{-1}=b, tbt^{-1}=(ab)^{-1}, zaz=a^{-1}, zbz=b^{-1}\rangle$$acting on $S$ with signature $(0; 2,2,3,3),$ 
\item isomorphic to  the group of order $6q^2$ presented as  $$H_2 =\langle a,b,t,z : a^q=b^q=t^3=z^2=[a,b]=1,$$ $$ tat^{-1}=b, tbt^{-1}=(ab)^{-1}, zaz=b, zbz=a, ztz=t^{-1}\rangle$$acting on $S$ with signature $(0; 2,2,3,3),$ or 
\item isomorphic to the group of order $12q^2$ presented as $$\hat{G}_1=\langle a,b,t,z,w : a^q=b^q=t^3=z^2=w^2=[a,b]=[z,w]=[z,t]=1,$$ $$ tat^{-1}=b, tbt^{-1}=(ab)^{-1}, zaz=a^{-1}, zbz=b^{-1}, waw=b, wtw=t^{-1} \rangle$$acting with  signature $(0; 2,2,2,3).$
\end{enumerate}

\end{theo}

\begin{theo} \label{t22} Let $q \geqslant 5$ be prime and set $g=1+q^2.$
The family $\mathscr{U}_g$ consists of at most $2q^2-q+1$ strata. Moreover, up to finitely many exceptions, the full automorphism group of $S \in \mathscr{U}_g$ is isomorphic to $G_2$ or isomorphic to the group of order $6q^2$ presented as$$\hat{G}_2=\langle a,t,z : a^{q^2}=t^3=z^2=[z,t]=1, tat^{-1}=a^s, zaz=a^{-1}\rangle$$where $s$ is a primitive third root of unity in $\mathbb{Z}_{q^2},$  acting with signature $(0; 2,2,3,3).$

\s
\s

\end{theo}

\begin{theo}\label{t23} Let $q \geqslant 5$ be prime  and set $g=1+q^2.$ The family $\mathscr{V}_g$ consists of only one stratum. Moreover,  the action of $G_3$ on each $S \in \mathscr{V}_g$ extends to an action of the group of order $6q^2$ presented as $$\hat{G}_3=\langle a, b,t, z : a^{q}=b^q=t^3=z^2=[a,b]=[z,t]=1,$$ $$\hspace{ 3 cm } tat^{-1}=a^{r}, tbt^{-1}=b^{r}, zaz=a^{-1}, zbz=b^{-1} \rangle,$$where $r$ is a primitive third root of unity in $\mathbb{Z}_q$ with signature $(0; 2,2,3,3).$ Furthermore, up to finitely many exceptions, the group $\hat{G}_3$ is isomorphic to the full automorphism group of $S.$
\end{theo}

We point out that 
the phrase {\it up to finitely many exceptions} in the theorem above cannot be deleted. In fact, to evidentiate that and for the sake of completeness, we shall construct, for each genus,  an explicit example of a Riemann surface lying in the family $\mathscr{U}_g$ and whose full automorphism group differs from the possibilities stated in the Theorem \ref{t22}.

\s
%
%

Let $S$ be a compact Riemann surface of genus $g \geqslant 2.$ We denote by $\mathscr{H}^1(S, \mathbb{C})$ the $g$-dimensional complex vector space of  1-forms on $S,$ and by $H_1(S, \mathbb{Z})$ the first integral homology group of $S.$
We recall that the Jacobian variety of $S,$ defined by $$JS=\mathscr{H}^1(S, \mathbb{C})^*/H_1(S, \mathbb{Z}),$$is an irreducible principally polarized abelian variety of dimension $g.$ The relevance of  the Jacobian variety  lies in  Torelli's  theorem, which establishes that two Riemann surfaces are isomorphic if and only if the corresponding Jacobian varieties are isomorphic as principally polarized abelian varieties. See, for example,  \cite[Section 11]{bl}.

\s

Let $H$ be a group of automorphisms of a compact Riemann surface $S$ and denote by $S_H$ the quotient compact Riemann surface given by the action of $H$ on $S.$ The  associated regular covering map $\pi : S \to S_H$   induces a homomorphism $$\pi^*: JS_H \to JS$$between the corresponding Jacobian varieties. The image $\pi^*(JS_H)$ is an abelian subvariety of $JS$ which is isogenous to $JS_H$. Thereby,    the well-known Poincar\'e's Reducibility theorem implies that there exists an abelian subvariety  of  $JS,$ henceforth  denoted by $\mbox{Prym}(S \to S_H)$ and called the Prym variety  associated to $\pi,$ such that  $$JS \sim JS_H \times \mbox{Prym}(S \to S_H),$$where $\sim$ stands for isogeny.

In what follows, we keep the same notations as in Theorem \ref{t1}.

\begin{theo} \label{t3} Let $q \geqslant 5$ and set $g=1+q^2.$

\s

If $S \in  \mathscr{C}_g$ then  the Jacobian variety $JS$ decomposes as $$JS \sim JS_{\langle a, b \rangle} \times JS_{\langle t \rangle}^3$$where $JS_{\langle a, b \rangle}$ is an abelian surface. Moreover, the Jacobian $JS_{\langle t \rangle}$ admits a further decomposition in terms of $\tfrac{q+1}{3}$ Prym varieties of the same dimension; concretely: 
\begin{equation} \label{arbol}JS_{\langle t \rangle} \sim \Pi_{n} \mbox{Prym}(S_{\langle ab^n \rangle} \to S_{\langle a, b \rangle})\end{equation}where $n$ runs over a subset of $\{1, \ldots, q-1\}$ which yields a maximal collection of pairwise non-conjugate subgroups of $G_1$ of the form $\langle ab^n \rangle$.

\s

If $S \in  \mathscr{U}_g$ then the Jacobian variety $JS$ decomposes as $$JS \sim JS_{\langle a \rangle} \times JS_{\langle t \rangle}^3$$where $JS_{\langle a \rangle}$ is an abelian surface. Moreover, $JS_{\langle t \rangle}^3$ admits a further decomposition in terms of a Prym variety of dimension $q(q-1)$ and the third power of a Jacobian variety of dimension $\tfrac{q-1}{3}$: $$JS_{\langle t \rangle}^3 \sim \mbox{Prym}(S \to S_{\langle a^q \rangle}) \times JS_{\langle a^q, t \rangle}^3.$$

\s

If $S \in  \mathscr{V}_g$ then the Jacobian variety $JS$ decomposes as $$JS \sim JS_{\langle a, b \rangle} \times \mbox{Prym}(S_{\langle a \rangle} \to S_{\langle a,b \rangle}) \times \Pi_{n=0}^{q-1} \mbox{Prym}(S_{\langle a^nb \rangle} \to S_{\langle a, b \rangle})$$where $JS_{\langle a, b \rangle}$ is an abelian surface. Moreover, each Prym variety above can be further decomposed as third power of a Jacobian. Concretely: $$ JS \sim JS_{\langle a, b \rangle} \times JS_{\langle a, t \rangle}^3 \times \Pi_{n=0}^{q-1}JS_{\langle a^nb, t \rangle}^3$$ where the dimension of  $JS_{\langle a, t \rangle}$ and of each $JS_{\langle a^nb, t \rangle}$ is $\tfrac{q-1}{3}.$
\end{theo}

\s

We recall that Recillas' trigonal construction ensures that the Jacobian variety of a tetragonal Riemann surface is isomorphic to the Prym variety of an unramified two-fold cover of a trigonal Riemann surface; see  \cite{recillas} and also \cite[Section 12.7]{bl}. Very recently in  \cite{clr}, Lange, Rodr\'iguez and the first author  somehow  generalized this fact by studying Riemann surfaces $Z$ with a group of automorphisms isomorphic to $$G=N \rtimes P \,\, \mbox{ where }\,\, N \cong C_2^{p-1} \, \mbox{and } \, P \cong C_p$$such that $Z \to Z_G \cong \mathbb{P}^1$ ramifies over only values marked with $p,$ where $p$ is prime. Concretely, under this conditions, if $T$ denotes the $2^{p-1}$-gonal Riemann surface $Z_P$ then they proved that the Jacobian $JT$ is isogenous to the product of $(2^{p-1}-1)/p$ Prym varieties of unbranched two-fold regular covers of the $p$-gonal Riemann surface $Z_N$.  In addition, it was proved that the dimension of the involved Prym varieties is the same and that the corresponding isogeny is induced by multiplication by $2^{p-2}.$ 

\s

We  recall that a Riemann surface $S$ lying in the family $\mathscr{C}_g$ admits the action of $$G_1 =\langle a,b,t \rangle \cong  C_q^2 \rtimes C_3 \,\,\, \mbox{ where } q  \mbox{ is prime,}$$and that Theorem \ref{t3} says that the Jacobian of the $q^2$-gonal Riemann surface $S_{\langle t \rangle}$ is isogenous to the product of $(q+1)/3$ Prym varieties (of the same dimension) of unbranched $q$-fold regular covers of the trigonal Riemann surface $S_{\langle a,b \rangle}$. Note that our case $q=2$ and the case $p=3$ in \cite{clr} agree and corresponds to the well-studied action of the alternating group of order 12 on genus $g=5.$

\s

Here, we slightly modify the arguments employed in \cite[Section 3]{clr} to prove the following theorem, which, in some sense, also constitutes a generalization of Recillas' trigonal construction.

\begin{theo}\label{isogenia} Let  $q \equiv -1 \mbox{ mod } 3$ be an odd prime number and set $g=1+q^2.$ If $S \in  \mathscr{C}_g$   then the isogeny$$\Pi_{n} \mbox{Prym}(S_{\langle ab^n \rangle} \to S_{\langle a, b \rangle}) \to JS_{\langle t \rangle}$$ of Theorem \ref{t3}  is induced by multiplication by $q$. In particular, its kernel is contained in the $q$-torsion points.
\end{theo}

This article is organized  as follows. In Section \ref{ss2} we shall briefly review the basic background: Fuchsian groups, group actions on Riemann surfaces, decomposition of Jacobian varieties and the equisymmetric stratification of the moduli space. The proof of the theorems will be given in Sections \ref{ss3}, \ref{ss4}, \ref{ss5} and \ref{ss6}. Finally, we include a couple of remarks in Section \ref{ss7}.

\s

{\bf Acknowledgments.} 
The second author wishes to thank Professor Herbert Lange for valuable conversations during his research stay at University of Erlangen.

\section{Preliminaries} \label{ss2}

\subsection{Group actions and Fuchsian groups} By a  Fuchsian group we mean a discrete group of automorphisms of the upper-half plane $$\mathbb{H}=\{z \in \mathbb{C}: \mbox{Im}(z) >0 \}.$$  

If $\Delta$ is a Fuchsian group and the orbit space $\mathbb{H}_{\Delta}$ given by the action of $\Delta$ on $\mathbb{H}$ is  compact, then the algebraic structure of $\Delta$ is determined by its so-called signature; namely, the tuple \begin{equation} \label{sig} \sigma(\Delta)=(h; m_1, \ldots, m_l),\end{equation}where $h$ is the genus of the quotient surface $\mathbb{H}_{\Delta}$ and $m_1, \ldots, m_l$ are the branch indices in the universal canonical projection $\mathbb{H} \to \mathbb{H}_{\Delta}.$ If $l=0$ then $\Delta$ is called a surface Fuchsian group.

\s

Let $\Delta$ be a Fuchsian group of signature \eqref{sig}. Then 
\begin{enumerate}
\item $\Delta$ has a canonical presentation, henceforth denoted by $\Delta(h; m_1, \ldots, m_l),$ given by generators $a_1, \ldots, a_{h}$, $b_1, \ldots, b_{h},$ $ x_1, \ldots , x_l$ and relations
\begin{equation}\label{prese}x_1^{m_1}=\cdots =x_l^{m_l}=\Pi_{i=1}^{h}[a_i, b_i] \Pi_{j=1}^l x_j=1,\end{equation}where $[u,v]$ stands for the commutator $uvu^{-1}v^{-1},$
\item the elements of $\Delta$ of finite order are conjugate to powers of $x_1, \ldots, x_l,$ and
\item the Teichm\"{u}ller space of $\Delta$ is a complex analytic manifold homeomorphic to the complex ball of dimension $3h-3+l$.
\end{enumerate}

Let $\Gamma$ be a group of automorphisms of $\mathbb{H}.$ If $\Delta$ is a subgroup of $\Gamma$ of finite index then $\Gamma$ is also Fuchsian and they  are related by the Riemann-Hurwitz formula $$2h-2 + \Sigma_{i=1}^l(1-\tfrac{1}{m_i})= [\Gamma : \Delta] (2h'-2 + \Sigma_{i=1}^s(1-\tfrac{1}{n_i})),$$where $\sigma(\Gamma)=(h'; n_1, \ldots, n_s).$ 

\s

Let $S$ be a compact Riemann surface and let $\mbox{Aut}(S)$ denote its full automorphism group. A  group $G$ is said to act on $S$ if there is a group monomorphism $G\to \Aut(S).$ The space of orbits $S_G$ of the action of $G$ on $S$ is endowed with a Riemann surface structure in such a way that the  projection $\pi_G:S \to S_G$ is holomorphic. 

\s

Compact Riemann surfaces and group actions can be understood in terms of Fuchsian groups as follows. By uniformization theorem, there is a surface Fuchsian group $\Gamma$ such that $S$ and $\mathbb{H}_{\Gamma}$ are isomorphic. Moreover, Riemann's existence theorem ensures that $G$ acts on $S \cong \mathbb{H}_{\Gamma}$ if and only if there is a Fuchsian group $\Delta$ containing $\Gamma$ together with a group  epimorphism \begin{equation*}\label{epi}\theta: \Delta \to G \, \, \mbox{ such that }  \, \, \mbox{ker}(\theta)=\Gamma.\end{equation*}

In this case, it is said that $G$ acts on $S$ with signature $\sigma(\Delta)$ and that this action is represented by the  surface-kernel epimorphism $\theta$. If $G$ is a subgroup of $G'$ then the action of $G$ on $S$ is said to extend to an action of $G'$ on $S$ if:\begin{enumerate}
\item there is a Fuchsian group $\Delta'$ containing $\Delta.$ 
\item the Teichm\"{u}ller spaces of $\Delta$ and $\Delta'$ have the same dimension, and
\item there exists a surface-kernel  epimorphism $$\Theta: \Delta' \to G' \, \, \mbox{ in such a way that }  \, \, \Theta|_{\Delta}=\theta.$$
\end{enumerate} 

An action is called maximal if it cannot be extended in the previous sense. A complete list of signatures of Fuchsian groups $\Delta$ and $\Delta'$ for which it may be possible to have an extension as before was determined by Singerman in \cite{singerman2}. 

\subsection{Actions and  stratification} \label{strati} Let $\text{Hom}^+(S)$ denote the group of orientation preserving homeomorphisms of $S.$ Two actions $\psi_i: G \to \mbox{Aut}(S)$  are termed topologically equivalent if there exist $\omega \in \Aut(G)$ and $f \in \text{Hom}^+(S)$ such that
\begin{equation}\label{equivalentactions}
\psi_2(g) = f \psi_1(\omega(g)) f^{-1} \hspace{0.5 cm} \mbox{for all} \,\, g\in G.
\end{equation}

Each homeomorphism $f$ satisfying \eqref{equivalentactions} yields an automorphism $f^*$ of $\Delta$ where $\mathbb{H}_{\Delta} \cong S_G$. If $\mathscr{B}$ is the subgroup of $\mbox{Aut}(\Delta)$ consisting of them, then $\mbox{Aut}(G) \times \mathscr{B}$ acts on the set of epimorphisms defining actions of $G$ on $S$ with signature $\sigma(\Delta)$ by $$((\omega, f^*), \theta) \mapsto \omega \circ \theta \circ (f^*)^{-1}.$$  

Two epimorphisms $\Delta \to G$ define topologically equivalent actions if and only if they belong to the same $(\mbox{Aut}(G) \times \mathscr{B})$-orbit (see \cite{Brou}; also \cite{Harvey} and \cite{McB}). 


\s

Let $\mathscr{M}_g$ denote the moduli space of compact Riemann surfaces of genus $g \geqslant2.$ It is well-known  that $\mathscr{M}_g$ is endowed with a structure of complex analytic space of dimension $3g-3,$ and that for $g \geqslant4$ its singular locus $\mbox{Sing}(\mathscr{M}_g)$ agrees with the set of points representing compact Riemann surfaces with non-trivial automorphisms. 

\s

According to Broughton \cite{b} (see also \cite{Harvey}), the singular locus $\mbox{Sing}(\mathscr{M}_g)$ admits an  equisymmetric stratification, where 
each equisymmetric stratum, if nonempty, corresponds to one topological class of maximal actions. More precisely:
\begin{enumerate}
\item the closure  $\bar{\mathscr{M}}_g^{G, \theta}$ of the equisymmetric stratum $\mathscr{M}_g^{G, \theta}$ consists of those Riemann surfaces of genus $g$ with an action of $G$ with fixed topological class given by $\theta,$
\item $\bar{\mathscr{M}}_g^{G, \theta}$ is a closed  irreducible algebraic subvariety of $\mathscr{M}_g,$  
 
\item if ${\mathscr{M}}_g^{G, \theta} \neq \emptyset$  then it is a smooth, connected, locally closed algebraic subvariety of $\mathscr{M}_g$ which is Zariski dense in $\bar{\mathscr{M}}_g^{G, \theta},$ 
\item there are finitely many distinct strata, and $$\mbox{Sing}(\mathscr{M}_g) = \cup_{G \neq 1, \theta} \bar{\mathscr{M}}_g^{G, \theta}.$$
\end{enumerate}

Let ${\mathcal{F}}$ be a family of compact Riemann surfaces of genus $g$ such that each of its members has a group of automorphisms isomorphic to $G.$  Then $${\mathcal{F}}=\cup_{ \theta} {\mathcal{F}}_g^{G, \theta}$$where the stratum  ${\mathcal{F}}_g^{G, \theta}$ consists of those members of the family admitting an 
action of $G$ with topological class $\theta.$

\subsection{Decomposition of Jacobians with group action} \label{jacos} Let $G$ be a finite group and  let $W_1, \ldots, W_r$ be its rational irreducible representations. For each $W_l$ we denote by $V_l$ a complex irreducible representation of $G$ associated to it. 

It is classically  known that if the group $G$ acts on a compact Riemann surface $S$ then this action  induces a $\mathbb{Q}$-algebra homomorphism $$\Phi : \mathbb{Q} [G] \to \mbox{End}_{\mathbb{Q}}(JS)=\mbox{End}(JS) \otimes_{\mathbb{Z}} \mathbb{Q},$$from the rational group algebra of $G$ to the rational endomorphism algebra of $JS.$

For each $ \alpha \in {\mathbb Q}[G]$ we define the abelian subvariety $$A_{\alpha} := {\textup Im} (\alpha)=\Phi (n\alpha)(JS) \subset JS$$where $n$ is some positive integer chosen in such a way that $n\alpha \in {\mathbb Z}[G]$.

The decomposition of $1 = e_1 + \cdots + e_r \in \mathbb{Q}[G]$, where each $e_l$ is a uniquely determined central idempotent (computed from $W_l$), yields an isogeny $$JS \sim A_{e_1} \times \cdots \times A_{e_r}$$
which is $G$-equivariant. See \cite{l-r}.

Additionally, there are idempotents $f_{l1},\dots, f_{ln_l}$ such that $e_l=f_{l1}+\dots +f_{ln_l}$  where  $n_l=d_{V_l}/s_{V_l}$ is the quotient of the degree $d_{V_l}$ of $V_l$ and its Schur index $s_{V_l}$.  These idempotents provide $n_l$ subvarieties of $JS$ which are isogenous between them; let $B_l$ be one of them, for every $l.$ Thus, we obtain the following isogeny
\begin{equation} \label{eq:gadec}
JS \sim_G B_{1}^{n_1} \times \cdots \times B_{r}^{n_r} 
\end{equation}
called the group algebra decomposition of $JS$ with respect to $G$. See \cite{cr}.

If the representations are labeled in such a way that $W_1(=V_1)$ denotes the trivial one (as we will do in this paper) then $n_1=1$ and $B_{1} \sim JS_G$.

\s

Let $H$ be a subgroup of $G$ and consider the associated regular covering map $\pi_H:S \to S_H.$ It was proved in \cite{cr} that the group algebra decomposition of $JS$ with respect to $G$ induces the following isogeny decomposition of $JS_H:$   \begin{equation*}
JS_H \sim  B_{1}^{{n}_1^H} \times \cdots \times B_{r}^{n_r^H} \,\,\, \mbox{ with } \,\,\, {n}_l^H=d_{V_l}^H/s_{V_l}
\end{equation*}where $d_{V_l}^H$ stands for the dimension of the vector subspace $V_l^H$ of $V_l$ consisting of those elements fixed under $H.$ 

The isogeny above provides a criterion to identify if a factor in the group algebra decomposition of $JS$ with respect to $G$ is isogenous to the Jacobian variety of a quotient of $S$ or isogenous to the Prym variety of an intermediate covering of $\pi_G$. More precisely, if two  subgroups $H \le H'$ of $G$ satisfy $$d_{V_i}^{H}- d_{V_i}^{H'}=s_{V_i}$$for some fixed $2 \le i \le r$ and $$d_{V_l}^{H} - d_{V_l}^{H'} = 0$$for all $l \neq i$
such that $\dim(B_{l}) \neq 0,$ then 
$$B_{i} \sim \mbox{Prym}(S_{H} \to S_{H'}).$$ Furthermore if, in addition, the genus of $S_{H'}$ is zero then $
B_{i} \sim JS_{H}.$ See also \cite{yo}.

\s

Assume that $(\gamma; m_1, \ldots, m_l)$ is the signature of the action of $G$ on $S$  and that this action is represented by the surface-kernel epimorphism $$\theta: \Delta(\gamma; m_1, \ldots, m_l) \to G,$$with $\Delta(\gamma; m_1, \ldots, m_l)$ as in \eqref{prese}. Following \cite[Theorem 5.12]{yoibero}, for $2 \le i \le r,$ the dimension of $B_{i}$ in \eqref{eq:gadec} is given by 
\begin{equation*}\label{dimensiones}
\dim (B_{i})=k_{V_i}\Big[d_{V_i}(\gamma -1)+\frac{1}{2}\Sigma_{k=1}^l (d_{V_i}-d_{V_i}^{\langle \theta(x_k) \rangle} )\Big]  \end{equation*} where $k_{V_i}$ is the degree of the extension $\mathbb{Q} \le L_{V_i}$ with $L_{V_i}$ denoting a minimal field of definition for $V_i.$

\s

The decomposition of Jacobian varieties with group actions has been extensively studied, going back to contributions of Wirtinger, Schottky and Jung; see \cite{SJ} and \cite{W}. For decompositions of Jacobians with respect to special groups, we refer to the articles \cite{Ba}, \cite{d1}, \cite{Do}, \cite{nos}, \cite{IJR}, \cite{families},  \cite{LR2}, \cite{PA}, \cite{d3}, \cite{kanirubiyo}, \cite{RCR} and \cite{Ri}.

\subsection*{\it Notation} We denote by $C_n$ the cyclic group of order $n.$

\section{Proof of theorem \ref{t1}} \label{ss3}Let $q \geqslant 7$ be a prime number. 
\subsection*{Existence} Assume that  $q \equiv 1 \mbox{ mod }3$ and choose $r$ and $s$ to be primitive third roots of unity in $\mathbb{Z}_q$ and $\mathbb{Z}_{q^2}$  respectively. Then the groups  $$G_2=\langle a, t : a^{q^2}=t^3=1, tat^{-1}=a^{s}\rangle $$$$G_3=\langle a, b, t : a^{q}=b^q=t^3=[a,b]=1, tat^{-1}=a^{r}, tbt^{-1}=b^r\rangle$$exist and do not depend on the choice of $r$ and $s$. Meanwhile, the group  
 $$G_1=\langle a,b, t: a^q=b^q=t^3=[a,b]=1, tat^{-1}=b, tbt^{-1}=(ab)^{-1}\rangle $$exists for each prime number $q.$ Note that the equality $$2((q^2+1)-1)=3q^2(-2+4(1-\tfrac{1}{3}))$$shows that the Riemann-Hurwitz formula is satisfied for a $3q^2$-fold covering map from a Riemann surface of genus $q^2+1$ onto the projective line with four branch values marked with 3. Thus, by Riemann's existence theorem, the existence of the families $\mathscr{C}_g,$ $\mathscr{U}_g$ and $\mathscr{V}_g$ follow directly after noticing that the correspondences sending $(x_1, x_2, x_3, x_4)$ to $$(t,t^2,at^2,b^{-1}t^2), \,\, (t,t,at^2, a^{-s}t^2) \,\, \mbox{ and } \,\,  (t,at,bt^2, a^{-1}b^{-r}t^2)$$define surface-kernel epimorphisms from the Fuchsian group $$\Delta(0;3,3,3,3)=\langle x_1, x_2, x_3, x_4 : x_1^3=x_2^3=x_3^3=x_4^3=x_1x_2x_3x_4=1\rangle$$ onto $G_1, G_2$ and $G_3$ respectively.

\subsection*{Signatures} Let $S$ be a compact Riemann surface of genus $1+q^2$ endowed with an action of a group $G$ of order $3q^2.$ We claim that, if $q \geqslant 11$ then the possible signatures for action of $G$ on $S$ are $(1; 3)$ and $(0; 3,3,3,3).$ If $q=7$ then, in addition to the previous ones, the signature can be $(0; 7,7,21).$ 

\s

Suppose the signature of the action of $G$ on $S$ to be $(\gamma; m_1, \ldots, m_l).$ Then the Riemann-Hurwitz formula says that \begin{equation} \label{rh}\tfrac{2}{3}=2\gamma-2+l-\Sigma_{i=1}^l\tfrac{1}{m_j}\end{equation}and therefore $\gamma=0$ or 1.  In the latter case it is straightforward to see that $l=1$ and $m_1=3.$ So, we assume $\gamma=0$ and therefore \eqref{rh} turns into $$\Sigma_{i=1}^l\tfrac{1}{m_j}=l-\tfrac{8}{3}.$$As each $m_j \geqslant 3$ we see that $l=3$ or $4.$ For $l=4$ we have $$\Sigma_{i=1}^4\tfrac{1}{m_j}=\tfrac{4}{3} \, \mbox{ and therefore the signature is }\, (0; 3,3,3,3).$$For $l=3$ we have \begin{equation} \label{azul}\Sigma_{i=1}^3\tfrac{1}{m_j}=\tfrac{1}{3} \, \mbox{ and therefore } m_j \geqslant q \,   \mbox{ for each } \, j.  \end{equation}

Then:
\begin{enumerate} 
\item if $q \geqslant 11$ then the  left-hand side of  \eqref{azul} is at most $3/11;$ a contradiction. 
\item if $q=7$ then the $m_1=m_2=7$ and $m_3=21.$
\end{enumerate}

\subsection*{Uniqueness.} Let $G$ be a non-abelian group of order $3q^2.$ By the classical Sylow's theorems, $G$  is isomorphic to a semidirect product $Q \rtimes C_3$ where $Q$ has order $q^2.$ Observe that if $Q$ is isomorphic to $C_{q^2}$ then $q \equiv 1 \mbox{ mod } 3$ and $G$ is isomorphic $G_2.$ 

We now assume that $Q \cong C_q^2$ and therefore $G$ admits a presentation $$G=\langle  a, b, t. : a^q=b^q=t^3=[a, b]=1, tat^{-1}=a^{n}b^m, tbt^{-1}=a^ub^v \rangle$$for suitable $0 \leqslant n,m,u,v < q.$ We have two cases to consider. 

\s

{\bf Case 1.} Assume that $tat^{-1}\in \langle a \rangle.$ If $tat^{-1}=a$ then, as $t$ has order three and as $[t,b] \neq 1,$ we notice that  $q \equiv 1 \mbox{ mod } 3$ and $tbt^{-1}=a^ub^r$ for some $0 \leqslant u < q.$ We can assume $u=0;$ in fact, otherwise, define $$\hat{a}:=a, \,  \,\hat{b}:=a^{u(r-1)^{-1}}b$$and note that $t\hat{a}t^{-1}=\hat{a}$ and $t\hat{b}t^{-1}=(\hat{b})^r.$ It follows that $G$ is isomorphic to $$H_2=\langle a, b, t : a^{q}=b^q=t^3= [a,b]=[a,t]=1,  tbt^{-1}=b^{r}\rangle$$provided that $q \equiv 1 \mbox{ mod }3.$

If $tat^{-1}=a^r$ then $q \equiv 1 \mbox{ mod }3$ and $tbt^{-1}$ equals either  $b^r,$ $a^ub,$ or $a^ub^{r^2}$  for some $1 \leqslant u < q.$ In the first case $G$ is isomorphic to $G_3.$ In the second case, if $$\hat{a}:=a^{u(1-r)^{-1}}b, \,  \,\hat{b}:=a$$ then $t\hat{a}t^{-1}=\hat{a}$ and $t\hat{b}t^{-1}=(\hat{b})^r;$ thus, $G$ is isomorphic to $H_2.$ In the last case define $$\hat{a}:=a, \,  \,\hat{b}:=a^{u(r^2-r)^{-1}}b$$and note that $t\hat{a}t^{-1}=(\hat{a})^r$ and $t\hat{b}t^{-1}=(\hat{b})^{r^2};$ thus, $G$ is isomorphic to $$H_1=\langle a, b, t : a^q=b^q=t^3=1, tat^{-1}=a^r, tbt^{-1}=b^{r^2} \rangle.$$

{\bf Case 2.} Assume $tat^{-1} \notin \langle a \rangle.$ Then we can assume $tat^{-1}=b$ and therefore  $tbt^{-1}=(ab)^{-1},$ or   $tbt^{-1}=a^{\epsilon}b^{-\epsilon^2}$ provided that $q \equiv 1 \mbox{ mod }3,$ where $\epsilon$ is a primitive third root of $-1$ in $\mathbb{Z}_q$. In the latter case write $$\hat{a}:=a^{(\epsilon^2-2\epsilon)/(1+\epsilon)}b \,\, \mbox{ and } \,\,\hat{b}:=a^{\epsilon}b^{(1-2\epsilon)/(1+\epsilon)}$$ and note that $t\hat{a}t^{-1}=\hat{b}$ and $t\hat{b}t^{-1}=(\hat{a}\hat{b})^{-1};$ thus, $G$ is isomorphic to $G_1.$

\s

In brief, if $G$ has order $3q^2$ and is non-abelian, then:
\begin{enumerate}
\item[(a)] if $q \equiv 1 \mbox{ mod } 3$ then $G$ is isomorphic to either $H_1, H_2, G_2$ or $G_3.$  
\item[(b)] if $q \equiv -1 \mbox{ mod } 3$ then $G$ is isomorphic to $G_1.$ 
\end{enumerate}
\s

Once the possible non-abelian abstract groups have been determined, we proceed to study each possible signature separately.

\s

{\bf Signature $(1;3)$.} Clearly, an abelian group cannot act with this signature. Besides, the commutator subgroup of each of the aforementioned non-abelian groups does not have elements of order three. Thus,  this signature  cannot be realized.

\s

{\bf Signature $(0;7,7,21)$ for $q=7$.} As $G_1, G_2, G_3$ and $H_2$ cannot be generated by two elements of order seven, we see that $G$ is necessarily abelian; then, isomorphic to  $C_3 \times C_{49}$ or $C_3 \times C_7^2.$ However, these  groups cannot be generated by two elements of order 7 with product of order 21. Thus,  this signature  cannot be realized.

\s

{\bf Signature $(0;3,3,3,3)$.} Note that $H_2$ and each  abelian group of order $3q^2$ cannot be generated by elements of order 3. On the other hand,  if $q \equiv 1 \mbox{ mod }3$ then the map  $$a \mapsto ab^{-r}, \, b \mapsto ab^{-r^2}$$defines a group isomorphism between $H_1$ and $G_1.$ Thus, if $q \equiv 1 \mbox{ mod } 3$ then $S$ belongs to $\mathscr{C}_g,$ $\mathscr{U}_g$ or $\mathscr{V}_g,$ and if $q \equiv -1 \mbox{ mod } 3$ then $S$ belongs to  $\mathscr{C}_g.$

\s

The proof of the theorem is done.

\section{Proof of Theorem \ref{t21}, \ref{t22} and \ref{t23}}\label{ss4}

For each surface-kernel epimorphism$$\theta: \Delta (0; 3,3,3,3) \to G_j \,\, \,\mbox{ for }\,\,\, j=1,2,3$$representing an action of $G_j$ on a compact Riemann surface $S$, we write$$g_i:=\theta(x_i) \,\,\, \mbox{ for each } \,\,\, i =1,2,3,4,$$and, for the sake of simplicity, we identify $\theta$ with the $4$-uple or   generating vector $$\theta=(g_1, g_2, g_3, g_4).$$ 

Note that the groups $G_1, G_2$ and $G_3$ have two conjugacy classes of elements of order three, represented by $t$ and $t^2.$ Then, as $g_1, g_2, g_3$ and $g_4$ have order three and as the product of them must be 1, we see that among them there are exactly two that are conjugate to $t.$ Moreover, after considering an inner automorphism of the group, we can assume $g_1=t,$ as we shall do in the sequel. 

\s

We record here that, by classical results on inclusions of Fuchsian groups due to Singerman (see \cite[Theorem 1]{singerman2}), the action of a group of order $d$ on a Riemann surface with signature $(0; 3,3,3,3)$ can be possibly extended to an action of a group of order $2d$ with signature $(0;2,2,3,3),$ and this action, in turn, can be possibly extended to  a maximal action of a group of order $4d$ with signature $(0;2,2,2,3)$.

We shall use repeatedly this fact in what follows.

\subsection*{Proof of Theorem \ref{t21}} Let $S \in \mathscr{C}_g.$ We recall that $S$ admits an action of $$G_1=\langle a, b, t : a^{q}=b^q=t^3=[a,b]=1, tat^{-1}=b, tbt^{-1}=(ab)^{-1}\rangle$$with signature $(0; 3,3,3,3).$ Assume that the action of $G_1$ on $S$ is given by  $$\theta=(t, a^ub^vt^2, g_3, g_4) \,\, \mbox{ for some }\,\, u,v \in \{0, \ldots, q-1\}$$and that $g_4$ is conjugate to $t.$ 
We have three cases to consider:

\s

{\bf Type 1:} $u =v=0.$ In this case $\theta$ is of the form $$(t, t^2, a^{n}b^{m}t^2, a^{m}b^{m-n}t) \,\, \mbox{ for some } \,\, m,n \in \{0, \ldots, q-1\}.$$Note that if $m=0$ then $n \neq 0$ and therefore, by sending $a$ and $b$ to an appropriate power of themselves, we obtain that $\theta$ is equivalent to  $$(t,t^2,at^2, b^{-1}t^2).$$Analogously,  if $m \neq 0$ then $\theta$ is equivalent to the surface-kernel epimorphism $$(t,t^2,a^nbt^2, ab^{1-n}t) \,\, \mbox{ for some }\,\, n \in \{0, \ldots, q-1\}.$$Thus, there are at most $q+1$ pairwise non-equivalent surface-kernel epimorphisms of type 1.
\s

{\bf Type 2:}  $u\neq 0$ and $v=0$ or $u=0$ and $v\neq 0.$  In this case, we can assume $u=1$ and $v=1$ respectively and therefore $\theta$ is equivalent to either  $$\mbox{type 2A}: \,\, (t,at^2, a^{n}b^{m}t^2, a^{m+1}b^{m-n+1}t) \,\, \mbox{ where }\,\, n,m \in \{0, \ldots, q-1\}.$$
or  $$\mbox{type 2B}: \,\, (t,bt^2, a^{n}b^{m}t^2, a^{m-1}b^{m-n}t) \,\, \mbox{ where }\,\, n,m \in \{0, \ldots, q-1\}.$$Thus,  there are at most $2q^2$ pairwise non-equivalent surface-kernel epimorphisms of type 2.

\s

{\bf Type 3:} $u, v \neq 0.$ In this case we can assume $v=1,$ and then $\theta$ is equivalent to $$(t, a^lbt^2, a^{n}b^{m}t^2, a^{l+m-1}b^{l-n+m}t)$$where $l \in \{1, \ldots, q\}$ and $n,m \in \{0, \ldots, q-1\}.$ Thus, there are at most $q^2(q-1)$ pairwise non-equivalent surface-kernel epimorphisms of type 3.

\s

All the above says that  $\mathscr{C}_q$ consists of at most $q^3+q^2+q+1$ strata.

\s

We now prove that, up to finitely many exceptions, the full automorphism group of $S$ is one of the possibilities stated in the theorem. 

\s

{\bf Claim 1.} There are actions of $G_1$ that extend to actions of $H_1, H_2$ and $\hat{G}_1$.
 
 \s

Consider the surface-kernel epimorphism $\Theta_i: \Delta(0; 2,2,3,3) \to H_i$ defined by $$\Theta_1(y_1)= az, \,\, \Theta_1(y_2) = zb, \,\, \Theta_1(y_3) = (ab)^{-1}t^2 \,\, \mbox{ and }\, \Theta_1(y_4) = t$$ $$\Theta_2(y_1)= a^{-1}bz, \,\, \Theta_2(y_2) = z, \,\, \Theta_2(y_3) = ab^{-1}t^{2} \,\, \mbox{ and }\, \Theta_2(y_4) = t.$$

The subgroup of $\Delta(0; 2,2,3,3)$ generated  by 
$$\hat{x}_1:=y_3, \,\, \hat{x}_2:=y_4, \,\, \hat{x}_3:=y_1y_3y_1 \,\, \mbox{ and } \,\, \hat{x}_4:=y_1y_4y_1$$is a Fuchsian group of signature $(0; 3,3,3,3).$ In addition $$\Theta_1(\hat{x}_1)=(ab)^{-1}t^2, \,\,\Theta_1(\hat{x}_2)=t, \,\, \Theta_1(\hat{x}_3)=a^3b^2t^2 \,\, \mbox{ and } \Theta_1(\hat{x}_4)=ab^{-1}t,$$$$\Theta_2(\hat{x}_1)=ab^{-1}t^2, \,\,\Theta_2(\hat{x}_2)=t, \,\, \Theta_2(\hat{x}_3)=a^{-1}b^4t \,\, \mbox{ and } \Theta_2(\hat{x}_4)=a^{-3}t^2,$$showing that  the restriction $$\Theta_i|_{\langle \hat{x}_1, \hat{x}_2, \hat{x}_3, \hat{x}_4\rangle}: \langle \hat{x}_1, \hat{x}_2, \hat{x}_3, \hat{x}_4\rangle \cong  \Delta(0; 3,3,3,3) \to \langle a, b,t \rangle \cong  G_1$$is equivalent to a surface-kernel epimorphism $\theta$ of type 2B for $i=1$ and of type 3 for $i=2$. It follows that there are two   equisymmetric strata whose members admit an action of $H_1$ and $H_2$ respectively, with signature $(0; 2,2,3,3).$ 
 
\s

Now, consider the surface-kernel epimorphism $\Theta: \Delta(0; 2,2,2,3) \to \hat{G}_1$ defined by $$\Theta(y_1)= az, \,\, \Theta(y_2) = zw, \,\, \Theta(y_3) = wt \,\, \mbox{ and }\, \Theta(y_4) = (at)^{-1}$$The subgroup of $\Delta(0; 2,2,2,3)$ generated by $$\hat{x}_1:=y_3y_4y_3, \,\, \hat{x}_2:=y_4, \,\, \hat{x}_3:=y_1y_3y_4y_3y_1 \,\, \mbox{ and } \,\, \hat{x}_4:=y_1y_4y_1$$is a Fuchsian group of signature $(0; 3,3,3,3).$ In addition $${\Theta}(\hat{x}_1)=a^{-1}t, \,\, {\Theta}(\hat{x}_2)=abt^2, \,\, {\Theta}(\hat{x}_3)=a^2b^{-1}t \,\, \mbox{ and } {\Theta}(\hat{x}_4)=at^{2}$$showing that  the restriction $$\Theta|_{\langle \hat{x}_1, \hat{x}_2, \hat{x}_3, \hat{x}_4\rangle}: \langle \hat{x}_1, \hat{x}_2, \hat{x}_3, \hat{x}_4\rangle \cong  \Delta(0; 3,3,3,3) \to \langle a,b,t \rangle \cong  G_1$$is equivalent to a surface-kernel epimorphism $\theta$ of type 3. It follows that there is an equisymmetric stratum whose members admit  an action of $\hat{G}_1$ with signature $(0; 2,2,2,3).$ 

This proves the claim.

\s

{\bf Claim 2.} Assume that an action of $G_1$ on $S$ extends to an action of a group $H$ of order $6q^2$ with signature $(0; 2,2,3,3).$ Then $H \cong H_1$ or $H \cong H_2.$

\s

Let $P=\langle t \rangle$ and let $K=\langle z \rangle$ be a Sylow $2$-subgroup of $H$. Observe that $P$ is a normal subgroup of $PK$ and therefore $ztz=t$ or $ztz=t^{-1}.$

\s

Assume $ztz=t.$ We write $zaz=a^{n}b^m$ and observe that the equalities $$t(zaz)t^{-1}=z(tat^{-1})z=zbz \,\, \mbox{ and }\,\, t(zaz)t^{-1}=t(a^nb^m)t^{-1}=b^n(ab)^{-m}$$show that $zbz=a^{-m}b^{n-m}.$ But,  as $z$ has order two, we obtain that $n=\pm 1$ and $m=0.$ It follows that either $$zaz=a, zbz=b \,\, \mbox{ or } zaz=a^{-1}, zbz=b^{-1}.$$
%

The former case  must be disregarded since the corresponding group has only one involution and therefore cannot act with signature $(0; 2,2,3,3)$. The latter one is  $H_1.$

\s

Assume $ztz=t^{-1}.$  We write $zaz=a^{n}b^m$ and --by arguing as in the previous case-- we notice that  $zbz=a^{m-n}b^{-n}.$ In addition, as $z$ has order two, we obtain that $$m^2-mn+n^2=1.$$Thus $(n,m)\in \{ (1,0),(1,1), (0,1),(0,-1), (-1,0), (-1,-1) \}.$ A routine computation shows that, for each case, it is possible to find elements $A$ and $B$ of order $q$ such that $$tAt^{-1}=B, \, tBt^{-1}=(AB)^{-1}, \, zAz=B \mbox{ and }  zBz=A.$$Hence, $H$ is isomorphic to $H_2$ and the proof of the claim is done.

\s
{\bf Claim 3.} Assume that the action of $G_1$ on $S$ extends to the action of a group $\hat{G}$ of order $12q^2$ with signature $(0; 2,2,2,3).$ Then $\hat{G} \cong \hat{G}_1.$
 
\s

Let $P=\langle t \rangle, Q=\langle a, b \rangle$ and $K$ be a Sylow $2$-subgroup of $\hat{G}$. Notice that $Q$ is a normal subgroup of $\hat{G}$. If $\hat{G}/Q$ is abelian or isomorphic to the alternating group then $QK$ is a normal subgroup of $\hat{G}$ and therefore the involutions of $\hat{G}$ are contained in $QK;$ this contradicts the fact that $\hat{G}$ can be generated by involutions. In addition, as $\hat{G}/Q$ is not isomorphic to the alternating group, we see that $G_1=QP$ is a normal subgroup of $\hat{G}.$ 

Besides, we notice that $K \cong C_2^2$ since otherwise if we write $K=\langle d \rangle$ then the normal subgroup $\langle d^2 \rangle G_1$ would contain all the involutions of $\hat{G}.$ Let us write $$K=\langle z,w : z^2=w^2=(zw)^2=1 \rangle.$$ Note that, as $P$ is normal in $PK$ and $PK$ is non-abelian, we have that $|PK/\mathbf{C}_{PK}(P)| \neq 2,$ where $\mathbf{C}_{PK}(P)$ denotes the centralizer of $P$ in $PK.$ Thus, we can write $$[z,t]=1 \,\, \mbox{ and }\,\, wtw=t^{-1}.$$

We now argue as done in Claim 2 to ensure that:

\begin{enumerate}
\item As $[z,t]=1$ we have that  $zaz=a^{-1}, zbz=b^{-1}$ or $zaz=a, zbz=b.$ 
\item As $wtw=t^{-1}$  we have that $waw=b$ and $wbw=a.$
\end{enumerate}

Finally, the latter possibility in (1) must be disregarded, since the subgroup $\langle a,b,t,z \rangle$ of order $6q^2$ cannot act with signature $(0; 2,2,3,3).$ The proof of the claim is done.

%
%
%
%
%
%
%

\subsection*{Proof of Theorem \ref{t22}} Let $S \in \mathscr{U}_g.$ We recall that $S$ admits an action of  $$G_2=\langle a, t : a^{q^2}=t^3=1, tat^{-1}=a^{s}\rangle$$where $s$ is a primitive third root of unity in $\mathbb{Z}_{q^2}$,  with signature $(0; 3,3,3,3).$ Assume that the action of $G_2$ on $S$ is given by  $$\theta=(t, a^vt, g_3, g_4)$$for some $v \in \{0, \ldots, q^2-1\}.$ 
We have three cases to consider:

\s

{\bf Type 1.} If $v=0$ then $\theta$ is of the form $$(t,t,a^nt^2, a^{-ns}t^2) \,\, \mbox{ where }1 \leqslant n < q^2$$with $n$ and $q^2$ coprime (otherwise $\theta$ is not surjective). Then, by passing to an appropriate power of $a$ we can assume $n=1$ and therefore $\theta$ is equivalent to $$\theta_1=(t,t,at^2, a^{-s}t^2).$$

\s

{\bf Type 2.} If $v=kq$ for some $k \in \{1, \ldots, q-1\}$ then $\theta$ is of the form $$(t,a^{kq}t,a^nt^2, a^{-ns-kq}t^2) \,\, \mbox{ where } 1 \leqslant n \ < q^2$$with  $n$ and $q^2$ are coprime.  Then, by passing to an appropriate power of $a$ we can assume $k=1$ and therefore $\theta$ is equivalent to $$\theta_{2, m}=(t,a^{q}t,a^{m}t^2, a^{-ms-q}t^2) \,\, \mbox{ for some } \,\, 1 \leqslant m < q^2$$ 
with $m$ and $q^2$ coprime; then, there are at most $q^2-q$ non-equivalent surface-kernel epimorphisms of this type. 
\s

{\bf Type 3.} If $v$ is different from 0 and is not a multiple of $q$ then we can assume $v=1$ and therefore $\theta$ is equivalent to the surface-kernel epimorphism $$\theta_{3,n}=(t,at,a^nt^2, a^{-1-ns}t^2) \,\, \mbox{ where } 0 \leqslant n < q^2$$Thus, there are at most $q^2$ non-equivalent surface-kernel epimorphisms of this type.

\s

All the above says that there are at most $2q^2-q+1$ strata.

\s

Consider the surface-kernel epimorphism $\Theta: \Delta(0; 2,2,3,3) \to \hat{G}_2$ defined by $$\Theta(y_1)= az, \,\, \Theta(y_2) = z, \,\, \Theta(y_3) = t \,\, \mbox{ and }\, \Theta(y_4) = a^{-s^2}t^2.$$The elements  
$$\hat{x}_1:=y_3, \,\, \hat{x}_2:=y_4, \,\, \hat{x}_3:=y_1y_3y_1 \,\, \mbox{ and } \,\, \hat{x}_4:=y_1y_4y_1$$generate a Fuchsian group of signature $(0; 3,3,3,3).$ In addition $${\Theta}(\hat{x}_1)=t, \,\, {\Theta}(\hat{x}_2)=a^{-s^2}t^2, \,\,{\Theta}(\hat{x}_3)=a^{1-s}t \,\, \mbox{ and } {\Theta}(\hat{x}_4)=at^{2}$$showing that  the restriction $$\Theta|_{\langle \hat{x}_1, \hat{x}_2, \hat{x}_3, \hat{x}_4\rangle}: \langle \hat{x}_1, \hat{x}_2, \hat{x}_3, \hat{x}_4\rangle \cong  \Delta(0; 3,3,3,3) \to \langle a,t \rangle \cong  G_2$$is equivalent to a surface-kernel epimorphism $\theta$ of type 3. It follows that there is an equisymmetric stratum whose members admit an action of $\hat{G}_2$ with signature $(0; 2,2,3,3).$

\s

{\bf Claim 1.} Assume that the action of $G_2$ on $S$ extends to the action of a group $\hat{G}$ of order $6q^2$ with signature $(0; 2,2,3,3).$ Then $\hat{G} \cong \hat{G}_2.$

\s

Observe that $\hat{G}$ is isomorphic to $$G_2 \rtimes \langle z : z^2=1 \rangle \,\, \mbox{ where } \,\, zaz=a^{\pm 1}, \,\, ztz=a^nt$$ for some $n.$ Note that if $zaz=a$ then $z$ is the unique involution and, consequently, the group cannot act with signature $(0; 2,2,3,3).$ Now, consider $t'=a^{n/2}t$ to see that $n$ can be chosen to be zero. The proof of the claim is done.
 
\s

{\bf Claim 2.} Each action of $\hat{G}_2$ on $S$ does not extend to any action of a group of order $12q^2$ with signature $(0;2,2,2,3).$

\s

To prove the claim we  proceed by contradiction; namely, we assume the existence of a group $H$ of order $12q^2$ such that:
\begin{enumerate}
\item $H$ contains a subgroup isomorphic to $\hat{G}_2$ and
\item $H$ is generated by three involutions.
\end{enumerate}

Observe that if $K$ is a Sylow $2$-subgroup of $H$ then $K \cong C_2^2.$ Indeed, as $\hat{G}_2$ has elements of order two and is normal in $H$, if $K$ were cyclic then the involutions of $H$ would be pairwise conjugate and this, in turn, would imply that $H$ cannot be generated by involutions. 

\s

On the other hand, as $Q=\langle a \rangle$ is the unique Sylow $q$-subgroup of $H,$ we see that $Q$ is normal in $H$ and the subset $QK$ turns into a subgroup of $H.$ Note that $QK$ is  non-normal in $H$ since otherwise all the involutions of $H$ would belong to $QK$.

\s

Let $P=\langle t \rangle$ and write $L=PK.$ Then, up to a conjugation, we can assume $H/Q \cong L.$  Note that $L$ has order 12 and that $P$ is normal in $L$ (otherwise, $K$ is normal in $L$ and therefore $QK$ would be normal in $H$). Without loss of generality, we can assume that $z \in K.$ We now choose an involution $w \in H$ in such a way  that $K=\langle z, w \rangle.$ Note that $[t,w] \neq 1,$ because $[t,z]=1$ and  $K$ is non-normal in $L.$ Moreover, as $P$ is normal in $L,$ we see that $$wtw=t^{-1}.$$

Finally, as $Q$ is normal in $H,$ we note that either $$waw^{-1}=a^{-1} \,\, \mbox{ or }\,\, waw^{-1}=a.$$As $tat^{-1}=a^{s}$  we see that:

\begin{enumerate}
\item in the first case $\mbox{Aut}(\langle a \rangle)$ contains a subgroup isomorphic to $\langle w ,t \rangle \cong \mathbf{S}_3;$ 
\item in the second case $(zw) a (zw)^{-1}=a^{-1}$ and therefore $\mbox{Aut}(\langle a \rangle)$ contains a subgroup isomorphic to $\langle zw, t \rangle \cong \mathbf{S}_3.$
\end{enumerate}The contradiction is obtained due to the fact that $\mbox{Aut}(\langle a \rangle)$ is abelian.

\subsection*{Proof of Theorem \ref{t23}} Let $S \in \mathscr{V}_g.$ We recall that $S$ admits an action of  $$G_3=\langle a, b, t : a^{q}=b^q=t^3=[a,b]=1, tat^{-1}=a^{r}, tbt^{-1}=b^r\rangle$$where $r$ is a primitive third root of unity in $\mathbb{Z}_q,$ acting with signature $(0; 3,3,3,3)$.  Assume that the action of $G_3$ on $S$ is given by the surface-kernel epimorphism $$\theta=(t, a^nb^mt, g_3, g_4) \,\, \mbox{ for some } \,\, n,m \in \{0, \ldots, q-1\} $$
We have two cases to consider:

\begin{enumerate}
\item[(a)] $n \neq 0$ and $m=0$ (this case is equivalent to the one in which $m \neq 0$ and $n=0$ by considering the automorphisms of $G$ given by $a \mapsto b, b \mapsto a$)
\item[(b)] $m$ and $n$ both different from zero (otherwise $\theta$ is not surjective)
\end{enumerate}

By sending $a$ and $b$ to appropriate power of themselves, we obtain that $\theta$ is equivalent to one of the following the surface-kernel epimorphism $$(t, at, g_3, g_4) \,\, \mbox{ or } \,\,  (t, abt, g_3, g_4).$$

These two epimorphisms  are, in turn, equivalent under the action of $a \mapsto ab, b \mapsto b.$ So, we can assume $\theta$ to be equivalent to$$(t, at, a^ub^v t^2, a^{-1-ru}b^{-rv}t^2)\,\, \mbox{ for some } 0 \leqslant u, v \leqslant q-1.$$ As $v \neq 0$ (otherwise $\theta$ is not surjective), again, by replacing $b$ by an appropriate power of it, we can assume $\theta$ to be equivalent to $$\theta_u:=(t, at, a^ub t^2, a^{-1-ru}b^{-r}t^2)\,\, \mbox{ for some } \,\, u \in \{0, \ldots, q-1\}.$$Note that the automorphism of $G$ given by $a \mapsto a, b \mapsto ab$ identifies $\theta_u$ with $\theta_{u+1};$ hence each $\theta$ is equivalent to the surface-kernel epimorphism $$\theta_0=(t, at, b t^2, a^{-1}b^{-r}t^2).$$

\s

{\bf Claim 1.} The action of $G_3$ extends to an action of $\hat{G}_3$ with signature $(0; 2,2,3,3).$ 

\s

Consider the surface-kernel epimorphism ${\Theta}: \Delta(0; 2,2,3,3) \to \hat{G}_3$ defined by $${\Theta}(y_1)=bz, \,\, {\Theta}(y_2)=a^rbz, \,\, {\Theta}(y_3)=t \,\, \mbox{ and } \,\, {\Theta}(y_4)=at^2$$and notice that $$\hat{x}_1:=y_3, \,\, \hat{x}_2:=y_4, \,\, \hat{x}_3:=y_1y_3y_1 \,\, \mbox{ and } \,\, \hat{x}_4:=y_1y_4y_1$$generate a Fuchsian group of signature $(0; 3,3,3,3).$ In addition $${\Theta}(\hat{x}_1)=t, \,\, {\Theta}(\hat{x}_2)=at^2, \,\, {\Theta}(\hat{x}_3)=b^{1-r}t \,\, \mbox{ and } {\Theta}(\hat{x}_4)=a^{-1}b^{1-r^2}t^2$$showing that the restriction of ${\Theta}$ to $\langle \hat{x}_1, \hat{x}_2, \hat{x}_3, \hat{x}_4\rangle$ defines a surface-kernel epimorphism which is, as the family consists of only one stratum, equivalent to  $\theta_0.$

\s

{\bf Claim 2.}  The action of $\hat{G}_3$ on $S$ does not extend to any action of a group of order $12q^2$ with signature $(0;2,2,2,3)$

\s

To prove the claim we proceed by contradiction; namely, we assume the existence of a group $H$ of order $12q^2$ such that:
\begin{enumerate}
\item $H$ contains a subgroup isomorphic to $\hat{G}_3=\langle a,b,t,z \rangle,$ and
\item $H$ is generated by three involutions.
\end{enumerate} 
 
We now argue similarly as done in the case of the family $\mathscr{U}_g.$ 
Let $K$ be a Sylow $2$-subgroup of $H$ and let $Q=\langle a, b \rangle.$ Then  $K \cong C_2^2$  and $QK$  is a non-normal subgroup of $H.$ If we write $P=\langle t \rangle$ and $L=PK$ then we can assume  that $H/Q \cong L$ and that $P$ is normal in $L.$ Up to  conjugation, we can suppose that $z \in K$ and that $K=\langle z, w \rangle$ for some involution $w \in H.$ In addition, as $P$ is normal in $L$ and $K$ is non-normal in $L,$ we see that $$wtw=t^{-1}.$$

As $Q$ is normal in $L$ and $w^2=1$, one of the following statements holds:$$[w,a]=1, \,\,\, waw=a^{-1}\,\,\, \mbox{ or } \,\,\,
 waw \in Q \setminus \langle a \rangle.$$

In the first case, we see that $$a^r=wa^rw=w(tat^{-1})w=(wtw)a(wt^{-1}w)=t^{-1}at=a^{r^2};$$a contradiction. In the second case we proceed analogously but considering $\hat{w}:=wz$ instead of $w;$ namely:$$a^r=\hat{w}a^r\hat{w}=\hat{w}(tat^{-1})\hat{w}=(\hat{w}t\hat{w})a(\hat{w}t^{-1}\hat{w})=t^{-1}at=a^{r^2}.$$

 Finally, if we assume that $waw \not\in  \langle a \rangle$ then $x:=a(waw)$ is non-trivial, belongs to $Q$ (and therefore $txt^{-1}=x^r$) and commutes with $w.$ Thus $$x^r=wx^rw=w(txt^{-1})w=t^{-1}(wxw)t=(wxw)^{r^2}=x^{r^2};$$a contradiction. The proof of the claim is done.

\section{Proof of Theorem \ref{t3}}\label{ss5}
 
Throughout this section we write $\omega_l=\mbox{exp}(\tfrac{2 \pi i}{l})$ for each $l \geqslant 2$ integer.

\subsection*{The family $\mathscr{C}_g$} Up to equivalence, the complex irreducible representations of $$G_1=\langle a, b, t : a^{q}=b^q=t^3=[a,b]=1, tat^{-1}=b, tbt^{-1}=(ab)^{-1}\rangle$$are: three of degree 1, given by $$\chi_k : a \mapsto 1, \,\,  b \mapsto 1, \,\, t \mapsto \omega_3^k \,\, \mbox{ for } \,\, k=0,1,2,$$and $\tfrac{q^2-1}{3}$ of degree $3,$ given by \begin{displaymath} \chi_{i,j} : a \mapsto \mbox{diag}(\omega_q^i, \omega_q^{j}, \omega_q^{-i-j}), \,\,  b \mapsto \mbox{diag}(\omega_q^j, \omega_q^{-i-j}, \omega_q^{i}),
\,\, t \mapsto
\left( \begin{smallmatrix}
0 & 1 & 0  \\
0 & 0 & 1 \\
1 & 0 & 0
\end{smallmatrix} \right)
\end{displaymath}with $$(i,j) \in (\{0, \ldots, q-1\}^2-\{(0,0)\})/R$$where $R$ is the relation given by  $(i,j) R (j, -i-j) R (-i-j,i).$

We  choose a set $\mathcal{P}$ of parameters $(i,j)$ in such a way that the set  $$\{\chi_{i,j}: (i,j) \in \mathcal{P}\}$$consists of a maximal collection of pairwise non-Galois associated representations of degree 3 of $G_1.$ Note that $G_1$ has $\tfrac{q+7}{3}$ conjugacy classes of cyclic subgroups; thus, the cardinality of $\mathcal{P}$ is $\tfrac{q+1}{3}.$

\s

Up to equivalence, the rational irreducible representations of $G_1$ are:
$$W_0 = \chi_0, \, W_1=\chi_1 \oplus \chi_2 \, \mbox{ and } \, W_{i,j}=\oplus_{\sigma}\chi_{i,j}^{\sigma}$$for each $(i,j) \in \mathcal{P},$ where the sum $\oplus_{\sigma}$ runs over the Galois group associated to the character field of $\chi_{i,j}$, which has order $q-1.$ 
\s

Let $S \in \mathscr{C}_g.$ Then the group algebra decomposition of $JS$ with respect to $G_1$ is\begin{equation} \label{pale}JS \sim B_{W_0} \times B_{W_1} \times \Pi_{(i,j) \in \mathcal{P}} B_{W_{i,j}}^3\end{equation}where $B_{W_0}=0.$ We apply \cite[Theorem 5.12]{yoibero} to ensure that, independently of the choice of the surface-kernel epimorphism representing the action of $G_1$ on $S$, the dimension of the abelian subvarieties in \eqref{pale} are
$$\dim(B_{W_1})=2 \,\, \mbox{ and } \,\, \dim (B_{W_{i,j}})=q-1$$for each $(i,j)	\in \mathcal{P}.$ Now, by \cite[Proposition 5.2]{cr} \begin{equation} \label{martes}JS_{\langle a, b \rangle} \sim B_{W_1} \,\, \mbox{ and } \,\, JS_{\langle t \rangle} \sim \Pi_{(i,j) \in \mathcal{P}}B_{W_{i,j}}\end{equation} showing that \eqref{pale} can be written as $$JS \sim JS_{\langle a, b \rangle} \times JS_{\langle t \rangle}^3.$$Note that $$\dim(JS_{\langle t \rangle})=\Sigma_{(i,j)\in \mathcal{P}} \dim(B_{W_{i,j}})=\tfrac{q^2-1}{3}.$$Now, if $$\mathcal{N}:=\{n  : in+j=0 \, \mbox{ for some }\, (i,j) \in \mathcal{P}\}$$then $\{\langle ab^n \rangle: n \in \mathcal{N}\}$ is a collection of maximal non-conjugate subgroups of $G_1$ of order $q.$  Clearly, the cardinality of $\mathcal{N}$ is $\tfrac{q+1}{3}$. Again, by \cite[Proposition 5.2]{cr}, for each $n \in \mathcal{N}$ we see that
 $$JS_{\langle ab^n \rangle} \sim B_{W_1} \times B_{W_{i_n, j_n}}$$where $(i_n, j_n)$ in the unique element of $\mathcal{P}$ such that  $i_n+nj_n=0.$   In particular  \begin{equation} \label{martes2}B_{W_{i,j}}\sim \mbox{Prym}(S_{\langle ab^n \rangle} \to S_{\langle a, b \rangle})\,\, \mbox{ where } \,\,  i+nj=0.\end{equation}
 
 Note that each covering map $S_{\langle ab^n \rangle} \to S_{\langle a, b \rangle}$ in unbranched and  $S_{\langle a, b \rangle}$ is trigonal. The desired isogeny follows from the second isogeny of \eqref{martes} together with \eqref{martes2}.

\subsection*{The family $\mathscr{U}_g$}
Up to equivalence, the complex irreducible representations of $$G_2=\langle a,  t : a^{q^2}=t^3=1, tat^{-1}=a^{s}\rangle,$$are: three of degree 1 given by $$ \chi_k : a \mapsto 1, \,\, t \mapsto w_3^k \,\, \mbox{ for } k=0,1,2,$$and $\tfrac{q^2-1}{3}$ of degree 3 given by  \begin{displaymath}V_{i}:  a \mapsto \mbox{diag}(\omega_{q^2}^i, \omega_{q^2}^{ir^2}, \omega_{q^2}^{ir}),
\,\, t \mapsto
\left( \begin{smallmatrix}
0 & 1 & 0  \\
0 & 0 & 1 \\
1 & 0 & 0
\end{smallmatrix} \right)
\end{displaymath}with $i \in \{1, \ldots, q^2-1\}/R$ where $R$ is the relation given by $i R ir R ir^2.$

Up to equivalence, the rational irreducible representations of $G_2$ are $$W_0=\chi_0, \,\, W_1 = \chi_1 \oplus \chi_2, \,\, W_2=\oplus_{\sigma}V_1^{\sigma} \,\, \mbox{ and } \,\, W_3=\oplus_{\sigma} V_q^{\sigma}$$where the sums $\oplus_{\sigma}$ are taken over the Galois group of the  character field of $V_1$ and $V_q,$ that have degree $\tfrac{q(q-1)}{3}$ and $\tfrac{q-1}{3}$ respectively. 

\s

If $S \in \mathscr{U}_g$ then the group algebra decomposition of $JS$ with respect to $G_2$ is  \begin{equation}\label{dia}JS \sim B_{W_0} \times B_{W_1} \times B_{W_{2}}^3 \times B_{W_{3}}^3.\end{equation}where $B_{W_0}=0.$ We apply \cite[Theorem 5.12]{yoibero} to ensure that, independently of the choice of the surface-kernel epimorphism representing the action of $G_2$ on $S$  $$\dim(B_{W_1})=2, \,  \dim(B_{W_{2}})=\tfrac{q(q-1)}{3} \, \mbox{ and }\, \dim(B_{W_{3}})=\tfrac{q-1}{3}$$
We now apply the results of \cite[Proposition 5.2]{cr} to obtain that \begin{equation}\label{verde}B_{W_1} \sim JS_{\langle a \rangle} \, \mbox{ and } \,  B_{W_2} \times B_{W_3} \sim JS_{\langle t \rangle}\end{equation}and the desired decomposition is obtained. Note that $$\dim(JS_{\langle t \rangle})=\tfrac{q(q-1)}{3}+\tfrac{q-1}{3}=\tfrac{q^2-1}{3}$$

Furthermore, we notice that $$B_{W_2}^3 \sim \mbox{Prym}(S \to S_{\langle a^q \rangle}) \, \mbox{ and } \, B_{W_3} \sim JS_{\langle a^q, t \rangle}$$and therefore, by the second isogeny of \eqref{verde}, we obtain  $$JS_{\langle t \rangle} \sim \mbox{Prym}(S \to S_{\langle a^q \rangle}) \times JS_{\langle a^q, t \rangle}^3.$$

\subsection*{The family $\mathscr{V}_g$} Up to equivalence, the complex irreducible representations of $$G_3=\langle a, b, t : a^{q}=b^q=t^3=[a,b]=1, tat^{-1}=a^{r}, tbt^{-1}=b^r\rangle$$ are: three of degree 1 given by $$ \chi_k : a \mapsto 1, \,\, b \mapsto 1, \,\, t \mapsto w_3^k \,\, \mbox{ for }\,\, k=0,1,2,$$and 
$\tfrac{q^2-1}{3}$ of degree 3 given by  \begin{displaymath}\chi_{i,j}:  a \mapsto \mbox{diag}(\omega_q^i, \omega_q^{ir^2}, \omega_q^{ir}), \,\,  b \mapsto \mbox{diag}(\omega_q^j, \omega_q^{jr^2}, \omega_q^{jr}),
\,\, t \mapsto
\left( \begin{smallmatrix}
0 & 1 & 0  \\
0 & 0 & 1 \\
1 & 0 & 0
\end{smallmatrix} \right)
\end{displaymath}with$$(i,j) \in \{0, \ldots, q-1\}^2-\{(0,0)\}/R$$where $R$ is the relation given by
 $(i,j)R(ir, jr)R(ir^2,jr^2).$
  
Up to equivalence, the rational irreducible representations of $G_3$ are $$W_0=\chi_0, \,\, W_1 = \chi_1 \oplus \chi_2, \,\, W_{{1,0}}=\oplus_{\sigma} \chi_{1,0}^{\sigma}$$ and $$W_{n,1}=\oplus_{\sigma} \chi_{n,1}^{\sigma} \,\, \mbox{ for each } \,\, 0\leqslant n \leqslant q-1$$where the sums $\oplus_{\sigma}$ are taken over the Galois group of the character field of $\chi_{1,0}$ and $\chi_{n,1}$ (of order $\tfrac{q-1}{3}$). Note that $G_3$ has $q+3$ conjugacy classes of cyclic subgroups.

\s

If $S \in \mathscr{V}_g$ then the group algebra decomposition of $JS$ with respect to $G_3$ is  \begin{equation} \label{h1}JS \sim B_{W_0} \times B_{W_1} \times B_{W_{{1,0}}}^3 \times \Pi_{n=0}^{q-1}B_{W_{n,1}}^3\end{equation} where $B_{W_0}=0.$  We now apply \cite[Theorem 5.12]{yoibero} to ensure that $$\dim(B_{W_1})=2 \,\, \mbox{ and } \,\, \dim(B_{W_{1,0}})=\dim(B_{W_{n,1}})=\tfrac{q-1}{3}$$for each $0 \leqslant n \leqslant q-1,$ and by \cite[Proposition 5.2]{cr} we deduce \begin{equation} \label{g1}JS_{\langle a,b \rangle} \sim B_{W_1}, \,\,\, JS_{\langle b \rangle} \sim B_{W_1} \times  B_{W_{1,0}}^3,  \,\,\, JS_{\langle a \rangle} \sim B_{W_1}  \times  B_{W_{0,1}}^3\end{equation}and $$ JS_{\langle a^nb \rangle} \sim B_{W_1}  \times  B_{W_{m,1}}^3$$where $m$ is such that $mn+1=0$ for each $1 \leqslant n \leqslant q-1.$ Hence \begin{equation} \label{h2}B_{W_{0,1}}^3 \sim \mbox{Prym}(S_{\langle a \rangle} \to S_{\langle a, b \rangle}) \,\ \mbox{ and } \,\,      B_{W_{1,0}}^3 \sim \mbox{Prym}(S_{\langle b \rangle} \to S_{\langle a, b \rangle}).\end{equation}and for each $n \in \{1, \ldots, q-1\}$ and $m$ as before we have that \begin{equation} \label{h3}B_{W_{n,1}}^3 \sim \mbox{Prym}(S_{\langle a^{m}b \rangle} \to S_{\langle a, b \rangle}).\end{equation} The first isogeny decomposition of $JS$ follows from the first isogeny in \eqref{g1} together with the isogenies \eqref{h1}, \eqref{h2} and \eqref{h3}. The last statement follows the fact that $$JS_{\langle a, t \rangle} \sim B_{W_{0,1}}, \,\, JS_{\langle b, t \rangle} \sim B_{W_{1,0}} \, \mbox{ and } \, JS_{\langle a^nb, t \rangle} \sim B_{W_{m,1}}$$where $1 \leqslant n \leqslant q-1$ and $m$ as before.

\section{Proof of Theorem \ref{isogenia}} \label{ss6}

Let $\phi: C \to C'$ be a covering map between compact Riemann surfaces. We recall that $\phi$ induces  two homomorphisms between the corresponding Jacobian varieties; namely, the norm  and the  pull-back $$N(\phi): JC \to JC' \,\, \mbox{ and } \,\, \psi^*: JC' \to JC.$$

For later use, we also keep in mind the fact that if $\phi$ is regular, say given by the action of a group of automorphisms $H$ of $C,$ then $$\phi^* \circ N(\phi) : JC \to JC \,\, \mbox{ is given by }\,\, z \mapsto \Sigma_{h \in H} h(z).$$

Let $S \in \mathscr{C}_g.$ We recall that $S$ admits the action of $$G_1=\langle a,b,t : a^q=b^q=t^3=[a,b]=1, tat^{-1}=b, tbt^{-1}=(ab)^{-1} \rangle.$$ 

Set $R:=S_{\langle a, b \rangle}$ and for each $n \in \{1, \ldots, q-1\}$  consider the   associated  maps $$f_n: S \to R_n:=S_{\langle ab^n\rangle}, \,\, \pi_n: R_n \to R \,\, \mbox{ and } \,\, \varphi : S \to T:=S_{\langle t \rangle}.$$

Thus, analogously as done in \cite{clr}, writing $$P_n:=\mbox{Prym}(\pi_n), \,\, A_n:=f_n^*(P_n), \,\, C_n:=N(\varphi)(A_n) \,\, \mbox{ and }\,\, B_n:=\varphi^*(C_n),$$we have that the following diagram \[
\xymatrix{
\,   & A_n \ar[rr]^{1+t+t^2} \ar[rd]^{N(\varphi)} & \, & B_n \ar[rr]^{\Sigma_{h \in \langle ab^n \rangle} h} \ar[rd]^{N(f_n)}& \, & A_n\\
P_n \ar[ru]^{f_n^*} \ar[rr]_{\alpha_i:=N(\varphi) \circ f_n^*}   & \, & C_n\ar[rr]_{ N(f_n)\circ \varphi^*} \ar[ru]^{\varphi^*} & \, & D_n \ar[ru]_{f_n^*}
}
\] is commutative, where $D_n:=N(f_n)(B_n).$

\s

We apply \cite[Proposition 2.2]{LR2} to the commutative diagram

\[
\xymatrix{
\,  & S \ar[rd]^{\varphi} \ar[ld]_{h:=\pi_n \circ f_n}& \, \\
R \ar[rd] & \, & T\ar[ld]\\
\,  & \mathbb{P}^1 & \,\\
}
\]to ensure that $\varphi^*(JT)$ is an abelian subvariety of $\mbox{Prym}(h).$ Then $$D_n=N(f_n)(\varphi^*(C_n)) \subset N(f_n)(\varphi^*(JT)) \subset N(f_n)(\mbox{Prym}(h)) \subset P_n,$$where the last inclusion follows from  $N(h)=N(\pi_n) \circ N(f_n).$ Note that \begin{equation}\label{key}\Phi_n \circ f_n^* = f_n^* \circ (N(f_n) \circ \varphi^* \circ N(\varphi) \circ f_n^*),\end{equation}where where $\Phi_n : A_n \to A_n$ is defined by $\Phi_n(z) = \Sigma_{h \in {\langle ab^n \rangle}}h \circ (1+t+t^2)(z).$

\s
{\bf Claim.} The map $$ N(f_n) \circ \varphi^* \circ N(\varphi) \circ f_n^*: P_n \to P_n$$ is the multiplication by $q.$

\s

To prove the claim we proceed analogously as done in the proof of \cite[Proposition 3.3]{clr}. By the equality \eqref{key} and since $f_n^*$ is an isogeny, we only need to verify that $\Phi_n$ is the multiplication by $q.$ Note that \begin{equation} \label{raton}\Phi_n(z)=\Sigma_{l=0}^{q-1}(ab^n)^l(z)+\Sigma_{l=0}^{q-1}(ab^n)^lt(z)+\Sigma_{l=0}^{q-1}(ab^n)^lt^2(z)\end{equation}for each $z \in JS.$ By \cite[Corollary 2.7]{rr2}, we can write $$A_n= \{ z \in JS : h(z)=z  \mbox{ for all } h \in {\langle ab^n \rangle} \mbox{ and } \Sigma_{j=0}^{q-1} b^j(z)=0\},$$and therefore if $z \in A_n$ then the first summand in  \eqref{raton} equals $qz.$ Besides$$(ab^n)^lt(z)= a^lb^{nl}t(z)=ta^{nl-l}b^{-l}(z);$$but $ab^n(z)=z$ implies that $(ab^n)^{nl-l}(z)=z$ and this, in turn, says that $a^{nl-l}b^{-l}(z)=b^{-l(n^2-n+1)}(z);$ consequently $$(ab^n)^lt(z)=tb^{s}(z) \,\, \mbox{ where } \,\, s:=-l(n^2-n+1).$$

Due to the fact that the polynomial  $x^2-x+1$ is irreducible provided that $q$ is not congruent to 1 modulo 3, one sees that $s$ runs over $\{0, \ldots, q-1\}$. Hence $$\Sigma_{l=0}^{q-1}(ab^n)^lt(z)=\Sigma_{s=0}^{q-1}tb^s(z)=t\Sigma_{s=0}^{q-1}b^s(z)=0,$$where the latter equality follows from the fact that $z \in A_n.$ Likewise,  the third summand in the right-hand side of \eqref{raton} equals zero, and the proof of the claim is done.

\s

Let $\mathcal{N}$ be a subset of $\{1, \ldots, q-1\}$ which yields a maximal collection of pairwise non-conjugate subgroups of $G_1$ the form $\langle ab^n \rangle$, and consider the isogeny $$\alpha:\Pi_{n \in \mathcal{N}} P_n \to JT$$of Theorem \ref{t3} which is induced by the addition map. 
If we denote by $\beta=(N(f_n))_{n \in \mathcal{N}}$ and $\alpha_n=N(\varphi) \circ f_n^*$ then $$\beta \circ \varphi^* \circ \alpha=\Pi_{n \in \mathcal{N}} (N(f_n^*) \circ \varphi^* \circ \alpha_n)$$and therefore, by the claim, the map $\beta \circ \varphi^* \circ \alpha$ is the multiplication by $q.$ It follows that the kernel of $\alpha$ is contained in the $q$-torsion points and the proof  is done.

\section{Some remarks} \label{ss7}

\subsection{A remark on the family $\mathscr{C}_g$}

If $q \equiv 1 \mbox{ mod } 3$ and $r$ is a primitive third root of unity in $\mathbb{Z}_q$ then the group $G_1$ of Theorem \ref{t1} is isomorphic to  $$\langle \alpha, \beta, \tau : \alpha^q=\beta^q=\tau^3=1, \tau \alpha \tau^{-1} = \alpha^r, \tau \beta \tau^{-1}=\beta^{r^2}\rangle$$(see the proof of Theorem \ref{t1}). 

By using the same arguments employed to prove Theorems \ref{t1} and \ref{t21}, the aforementioned presentation of $G_1$ permits us to prove that, if $q \equiv 1 \mbox{ mod } 3$, then:
\begin{enumerate}
\item the family $\mathscr{C}_g$ actually consists of at most $q^2+2q+1$ strata. 
\item if $S \in \mathscr{C}_g$ then $$JS \sim JS_{\langle \alpha, \beta \rangle} \times JS_{\langle \tau \rangle}^3$$where $JS_{\langle \alpha, \beta \rangle}$ is an abelian surface and $JS_{\langle \tau \rangle}$ decomposes in terms of $\tfrac{q-1}{3}$ Prym varieties of the same dimension $q-1$ and two Jacobians of dimension $\tfrac{q-1}{3}.$ Concretely  
\begin{equation*} JS_{\langle \tau \rangle} \sim JS_{\langle \alpha, \tau \rangle} \times JS_{\langle \beta, \tau \rangle} \times \Pi_{n} \mbox{Prym}(S_{\langle \alpha \beta^n \rangle} \to S_{\langle \alpha, \beta \rangle})\end{equation*}where $n$ runs over a subset of $\{1, \ldots, q-1\}$ which yields a maximal collection of pairwise non-conjugate subgroups of the form $\langle \alpha \beta^n \rangle$. 
\end{enumerate}

\subsection{On the phrase {\it up to finitely many exceptions} in Theorem \ref{t21}, \ref{t22} and \ref{t23}} 

Let $$\Delta(0; 2,6,6)=\langle y_1, y_2, y_3 : y_1^2=y_2^6=y_3^6=y_1y_2y_3=1 \rangle$$be a triangle Fuchsian group of signature $(0; 2,6,6),$ and consider the group  $$\tilde{G}_2=\langle a,c, \omega : a^{q^2}=c^6=\omega^2=1, cac^{-1}=a^{\epsilon}, [\omega, a]=[\omega, c]=1\rangle$$where $\epsilon$ is a primitive sixth root  of unity in  $\mathbb{Z}_{q^2}.$
\s

Note that $\tilde{G}_2$ has order $12q^2$ and the rule $\eta: \Delta(0; 2,6,6) \mapsto \tilde{G}_2$ given by $$(y_1, y_2, y_3) \mapsto (\omega c^3, a^{-\epsilon^2}c^2\omega, ac)$$is a surface-kernel epimorphism of type $(0; 2,6,6).$ It follows that $\eta$ guarantees the existence of a Riemann surface $X$ of genus $1+q^2$ with an action of $\tilde{G}_2$ given by $\eta.$

\s

To see that $X$ belongs to $\mathscr{U}_g$ it is enough to consider the action on $X$ of the subgroup of $\tilde{G}_2$ given by $\langle a, b:=c^2 \rangle \cong G_2.$  In conclusion, $X \in \mathscr{U}_g$ and $X$ has strictly more automorphisms than $6q^2.$

\end{document}